\documentstyle[twocolumn]{article}
\def\dsup{\displaystyle\sup}

\def\dinf{\displaystyle\inf}

\makeatletter
\def\@normalsize{\@setsize\normalsize{12pt}\xpt\@xpt
\abovedisplayskip 10pt plus2pt minus5pt\belowdisplayskip \abovedisplayskip
\abovedisplayshortskip \z@ plus3pt\belowdisplayshortskip 6pt plus3pt
minus3pt\let\@listi\@listI}

\def\subsize{\@setsize\subsize{12pt}\xipt\@xipt}

\def\section{\@startsection {section}{1}{\z@}{24pt plus 2pt minus 2pt}
{12pt plus 2pt minus 2pt}{\large\bf}}

\def\subsection{\@startsection {subsection}{2}{\z@}{12pt plus 2pt minus 2pt}
{12pt plus 2pt minus 2pt}{\subsize\bf}}
\makeatother

\begin{document}
\date{}
\title{\bf {\sc Robust} SPR {\sc Synthesis for Low-Order Polynomial Segments and
 Interval Polynomials}
}
\author{{\bf Long Wang$^{1)}$ and
Wensheng Yu$^{2)}$} \\
\\
\it \small $^{1)}$ Dept. of Mechanics and Engineering Science, Center for Systems and Control,\\
\it \small Peking University, Beijing 100871, China, E-mail: longwang@mech.pku.edu.cn \\
 \it \small $^{2)}$  Lab of Engineering Science for Complex Systems, Institute of Automation, \\
  \it \small Chinese Academy of Sciences, Beijing 100080, China, E-mail: yws@compsys.ia.ac.cn \\
} \maketitle \thispagestyle{empty} \vspace{-0.5cm} {\em \small
\noindent {\bf Abstract: } We prove that, for low-order ($n\leq
4$) stable polynomial segments or interval polynomials, there
always exists a fixed polynomial such that their ratio is
SPR-invariant, thereby providing a rigorous proof of Anderson's
claim on SPR synthesis for the fourth-order stable interval
polynomials. Moreover, the relationship between SPR synthesis for
low-order polynomial segments and SPR synthesis for low-order
interval polynomials is also discussed.

\noindent {\bf Keywords:} Synthesis Method, Strict Positive
Realness(SPR), Constructive Design, Robust Stability,
SPR-invariance, Transfer Functions, Polynomial Segments, Interval
Polynomials, Polytopic Polynomials.}
\begin{center}
\section*{1. INTRODUCTION}
\end{center}
The notion of strict positive realness (SPR) of transfer functions plays an important role
in absolute stability theory, adaptive control and system identification[1-5].
In recent years, stimulated by the robustness analysis method[6-9],
the study of robust strictly positive real transfer functions
has received much attention, and great progress has been made[10-23].
However, most available results belong to the category of robust SPR analysis. Much work remains to be
done in  robust SPR synthesis.

Synthesis problems are mathematically more difficult than analysis problems.
Usually, the synthesis problems require answering questions of existence and construction,
whereas the analysis problems can be dealt
with under the assumption of existence.
Synthesis problems are of more practical significance
from the engineering application viewpoint.

The basic statement of the robust strictly positive real synthesis problem is as follows:
Given an $n$-th order robustly stable polynomial set $F$, does there
exist, and how to
construct a (fixed) polynomial $b(s)$ such that, $\forall a(s)\in F,$  $a(s)/b(s)$ is strictly positive
real? (If such a polynomial $b(s)$ exists,
then we say that $F$ is synthesizable.)

When $F$ is a low-order ($n\leq 3$) interval polynomial set,
the synthesis problem above
has been considered by a number of authors and several
important results$^{[13,14,16,17,19-21]}$ have been presented.
But when $F$ is a high-order ($n\geq 4$) interval
polynomial set, even in the case of $n=4, $
the synthesis problem above is still open$^{[16,17,19-21]}$.

By the definition of SPR, it is easy to know that the Hurwitz stability of $F$
is a necessary condition for the existence of polynomial $b(s)$.
In [13-15], it was proved that, if all polynomials in $F$
have the same even (or odd) parts, such a polynomial $b(s)$ always exists;
In [13,14,16,19-21],
it was proved that, if $n\leq 3$ and $F$ is a stable interval polynomial set,
such a polynomial $b(s)$ always exists;
Recent results in [18-20] show that, if $n\leq 3$ and $F$ is
the stable convex combination of two polynomials $a_1(s)$ and $a_2(s)$,
such a polynomial $b(s)$ always exists.
Some sufficient condition for robust SPR synthesis
are presented in [10,17,19-21],
especially, the design method in [19,20] is numerically
efficient for high-order polynomial segments and
interval polynomials, and the derived conditions in [19,20]
are necessary and sufficient for robust SPR synthesis
of low-order
($n\leq 3$) polynomial segments or interval polynomials.

It should be pointed out that, Anderson et al. [16] transformed the robust SPR synthesis
problem for the fourth-order interval polynomial set into linear programming problem
in 1990 (namely, equations (58)-(60) in [16]), and by using linear programming techniques,
they concluded
that such a  linear programming problem always had a solution, thus, it was thought that
the robust SPR synthesis problem for the fourth-order interval polynomial set had been solved.
But in 1993, a synthesizable example in [17] showed that the corresponding
linear programming problem  had no solution. Hence, for the fourth-order interval
polynomial set, on one hand, we could not prove theoretically the existence of robust SPR synthesis,
on the other hand, we could not find a counterexample that is not synthesizable. Therefore,
the robust SPR synthesis problem for interval polynomial set, even in the case of $n=4$, is still an open problem$^{[16,17,13,14,19-21]}$.

In this paper, we prove that, for low-order ($n\leq 4$) stable
polynomial segments or interval polynomials, there always exists a
fixed polynomial such that their ratio is SPR-invariant, thereby
providing a rigorous proof of Anderson's claim on SPR synthesis
for the fourth-order stable interval polynomials. Moreover, the
relationship between SPR synthesis for low-order polynomial
segments and SPR synthesis for low-order interval polynomials is
also discussed. Our proof is constructive, and is useful in
solving the general SPR synthesis problem. \vskip 6mm
\begin{center}
\section*{2. MAIN RESULTS}
\end{center}
\vskip -7mm
In this paper, $ P^n$ stands for the set of $n$-th
order polynomials with real coefficients, $ R$ stands for the
field of real numbers, $ \partial (p)$ stands for the order of
polynomial $ p(\cdot )$, and $H^n\subset P^n$ stands for the set
of $n$-th order Hurwitz stable polynomials.

In the sequel, $ p(\cdot )\in P^m, \,\, q(\cdot )\in P^n,\,\,
f(s)=p(s)/q(s)$ is a rational function.

{\bf  Definition 1}$^{[10,16,17,23]}$ \ \
A biproper rational function $ f(s)$ (i.e., $ \partial (p)=\partial (q)$) is said to be strictly positive real(SPR), if

(i) $f(s)$ is analytic in $ \mbox{Re} [s]\geq 0$, i.e., $ q(\cdot )
\in H^n$;

(ii) $ \mbox{Re}[f(j \omega )] > 0,\ \ \ \forall \omega \in R.$

If $f(s)=p(s)/q(s)$ is proper, it is easy to get the following property:

{\bf Lemma 1}$^{[11]}$\ \
If $f(s)=p(s)/q(s)$ is a proper rational function, $q(s)\in H^n,$ and
$\forall \omega \in R, \mbox{Re}  [f(j\omega )] > 0,$ then
$p(s)\in H^n\cup H^{n-1}.$

Denote $F=\{a_i(s)=s^n+\sum_{l=1}^n a_l^{(i)}s^{n-l},i=1,2 \} $
as the two endpoint polynomials of a stable polynomial segment $\overline F$ (convex combination),
it is easy to prove that:

{\bf Lemma {\bf 2}}$^{[16]}$\ \
$\forall a(s)\in \overline F,$  $b(s)/a(s)$ is strictly positive real, if and only if,
$b(s)/a_i(s),i=1,2,$ are strictly positive real.

Consider an interval polynomials
$$K=\{a(s)= s^n+\sum_{i=1}^na_is^{n-i},a_i\in [a_i^{-},a_i^{+}],i=1,2,\cdots,n\}
$$
Denote $F=\{a_i(s)=s^n+\sum_{l=1}^na_l^{(i)}s^{n-l},i=1,2,3,4 \} $
as the four Kharitonov vertex polynomials of $K$$^{[6-9]}$.

{\bf Lemma {\bf 3}}$^{[6]}$\ \
$K$ is robustly stable if and only if $a_i(s)\in H^n,i=1,2,3,4.$

The following result was proved by Dasgupta and Bhagwat$^{[10]}$:

{\bf Lemma {\bf 4}}$^{[10]}$\ \
$\forall a(s)\in K,$  $b(s)/a(s)$ is strictly positive real, if and only if,
$b(s)/a_i(s),i=1,2,3,4,$ are strictly positive real.

First, for a low-order ($n\leq 3$) stable convex combination of
polynomials, by [18-20], we have

{\bf Theorem 1}$^{[18]}$\ \
If $F=\{a_i(s)=s^n+\sum_{l=1}^na_l^{(i)}s^{n-l},$
$i=1,2.\} $
is the set of the two endpoint polynomials of a low order
($n\leq 3$) stable segment of polynomials (convex combination) $\overline{F}$,
then there always exists a fixed polynomial $b(s)$
such that $\forall a(s)\in \overline{F},$  $b(s)/a(s)$ is strictly positive real.

Furthermore, if $F$ is the four Kharitonov vertex polynomials of
a low-order ($n\leq 3$) stable interval polynomial set, then we have

{\bf Theorem 2}$^{[13,14,16,17,19-21]}$\ \
If $F=\{a_i(s)=s^n+\sum_{l=1}^na_l^{(i)}s^{n-l},$
$i=1,2,3,4.\} $
is the set of the four Kharitonov vertex polynomials of a low order
($n\leq 3$) stable
interval polynomial family $K$,
then there always exists a fixed polynomial $b(s)$
such that $\forall a(s)\in K,$  $b(s)/a(s)$ is strictly positive real.

The following two theorems are the main results of this paper:

{\bf Theorem 3}\ \
If $F=\{a(s)=s^4+a_1s^3+a_2s^2+a_3s+a_4,$
$b(s)=s^4+b_1s^3+b_2s^2+b_3s+b_4\} $
is the set of the two endpoint polynomials of a fourth order
stable segment of polynomials (convex combination),
then there always exists a fixed polynomial $c(s)$
such that $c(s)/a(s)$ and $c(s)/b(s)$ are strictly positive real.

Consider the fourth-order interval polynomials
$$
K=\left \{
\begin{array}{r}
a(s)= s^4+a_1s^3+a_2s^2+a_3s+a_4,\\
 a_i\in [a_i^{-},a_i^{+}], i=1,2,3,4
 \end{array}
 \right \}$$
Denote
$$
\begin{array}{l}
  a_1(s)= s^4+a_1^{+}s^3+a_2^{+}s^2+a_3^{-}s+a_4^{-} \\
  a_2(s)= s^4+a_1^{-}s^3+a_2^{-}s^2+a_3^{+}s+a_4^{+} \\
  a_3(s)= s^4+a_1^{+}s^3+a_2^{-}s^2+a_3^{-}s+a_4^{+} \\
  a_4(s)= s^4+a_1^{-}s^3+a_2^{+}s^2+a_3^{+}s+a_4^{-}
  \end{array}
  $$
as the four Kharitonov vertex polynomials of $K$$^{[6-9]}$.

{\bf Theorem 4}\ \
If $F=\{a_i(s),i=1,2,3,4\} $
is the set of the four Kharitonov vertex polynomials of a fourth order
 stable interval polynomial family,
then there always exists a fixed polynomial $b(s)$
such that $\forall a(s)\in F,$  $b(s)/a(s)$ is strictly positive real.

Note that in Theorem 3, $c(s)/a(s)$ and $c(s)/b(s)$ being strictly positive real implies
$\forall \lambda \in [0, 1], $ $\frac{c(s)}{\lambda a(s) + (1- \lambda) b(s)}$ being strictly positive real (by Lemma 2);
similarly, in Theorem 4, $\forall a(s)\in F,$  $b(s)/a(s)$ being strictly positive real
implies $\forall a(s)\in K,$  $b(s)/a(s)$ being strictly positive real (by Lemma 4).
\begin{center}
\section*{3. PROOFS OF MAIN RESULTS}
\end{center}
In order to prove the main results above, we must first establish some lemmas.

{\bf Lemma 5}\ \
Suppose $ a(s)=s^4+a_1s^3+a_2s^2+a_3s+a_4\in H^4,$
then the following quadratic curve is an  ellipse
in  the first quadrant of the $x$-$y$ plane:
$$(a_2^2-4a_4)x^2+2(2a_3-a_1a_2)xy+a_1^2y^2$$
$$\quad\quad\quad -2(a_2a_3-2a_1a_4)x-2a_1a_3y+a_3^2=0 $$
and this ellipse is tangent with $y$ axis at
$(0,\displaystyle\frac {a_3}{a_1}),$
tangent with the lines $x=a_1$ and $a_3y-a_4x=0$
at $(a_1,a_2-\displaystyle\frac {a_3}{a_1})$
and $(\displaystyle\frac {a_3^2}{a_2a_3-a_1a_4},
\displaystyle\frac {a_3a_4}{a_2a_3-a_1a_4})$, respectively.

{\bf Proof:}\ \ Since $ a(s)$ is Hurwitz stable, Lemma 5
can be verified by a direct calculation.

Let $a(s)=s^4+a_1s^3+a_2s^2+a_3s+a_4\in H^4,$ for notational
simplicity, denote
$$
\begin{array}{ll}
\Omega_e^a:=& \;\{(x,y)|(a_2^2-4a_4)x^2+2(2a_3-a_1a_2)xy+a_1^2y^2\\
& \;\quad\quad\quad
-2(a_2a_3-2a_1a_4)x-2a_1a_3y+a_3^2 <0\}\\
\Omega
_t^a:=& \;\{(x,y)| a_1-x \geq 0,a_2x-a_1y-a_3\geq 0,\\
& \;\quad\quad\quad\quad a_3y-a_4x > 0\}\\
\Omega ^a :=& \;\Omega _e^a \cup \Omega _t^a \\
\end{array}
$$
apparently, $\Omega ^a $ is a bounded convex set
in the $x$-$y$ plane.

{\bf Lemma 6}\ \
Suppose $ a(s)=s^4+a_1s^3+a_2s^2+a_3s+a_4\in H^4$
and $(x,y)\in \Omega ^a , $
let $ c(s):=s^3+xs^2+ys+\varepsilon,$ where $\varepsilon $
 is positive and sufficiently small,then
$\forall \omega \in R,\mbox{Re}
[\displaystyle\frac {c(j\omega )}{a(j\omega )}]>0.$

{\bf Proof:}\ \
Suppose $(x,y)\in \Omega ^a , $
let $ c(s):=s^3+xs^2+ys+\varepsilon,$ where $
\varepsilon >0 $ and is sufficiently small.

$ \forall \omega \in R,$ consider
$$
\begin{array}{ll}
\mbox{Re}
[\displaystyle\frac {c(j\omega )}{a(j\omega )}]
& \;=\displaystyle\frac  1{\mid
a(j\omega )\mid ^2}[(a_1-x)\omega ^6\\
& \;+(a_2x-a_1y-a_3)\omega ^4+(a_3y-a_4x)
\omega ^2\\
& \;+\varepsilon(\omega ^4-a_2\omega ^2+a_4)]\\
\end{array}
$$
In order to prove that $\forall \omega \in R,
\mbox{Re}
[\displaystyle\frac {c(j\omega )}{a(j\omega )}]>0,$
let $t= \omega ^2,$ we only need to prove that, for any sufficiently small
$\varepsilon >0$,
$$
\begin{array}{ll}
f(t):=& \;t[(a_1-x)t ^2+(a_2x-a_1y-a_3)t+(a_3y-a_4x)]\\
& \;\quad +
\varepsilon(t^2-a_2t+a_4)>0,\forall t\in [0,+\infty).\\
\end{array}$$
Since $(x,y)\in \Omega ^a,$ by definition of $\Omega ^a $ and Lemma 5,
$(x,y)$ satisfies $a_1-x>0,a_3y-a_4x>0,$ and
\begin{center}
$[a_2x-a_1y-a_3]^2-4(a_1-x)(a_3y-a_4x)< 0$
\end{center}
or
$$ a_1-x \geq 0,a_2x-a_1y-a_3\geq 0,
a_3y-a_4x> 0$$
\noindent Thus, $\forall t\in [0,+\infty)$
$$(a_1-x)t ^2+(a_2x-a_1y-a_3)t+(a_3y-a_4x)>0.$$

On the other hand, we have $f(0)>0,$ and for any $\varepsilon >0,$
if $t$ is a sufficiently large or sufficiently small positive number, we have $f(t)>0.$
Namely, there exist
$ 0<t_1<t_2$ such that, for all $ \varepsilon >0,
t\in [0,t_1]\cup [t_2,+\infty)$, $f(t)>0.$

Denote
$$ M=\dinf_{t\in[t_1,t_2]}t[(a_1-x)t ^2+(a_2x-a_1y-a_3)t+(a_3y-a_4x)],$$
$$ \quad N=\dsup_{ t\in[t_1,t_2]}{\mid t^2-a_2t+a_4 \mid }$$
then $M>0 $ and $N>0.$ Choosing
$0<\varepsilon <\displaystyle\frac {M}{N},$ by a direct calculation, we have
$$
\begin{array}{ll}
f(t)=& \;t[(a_1-x)t ^2+(a_2x-a_1y-a_3)t+(a_3y-a_4x)]\\
& \;\quad +\varepsilon(t^2-a_2t+a_4)>0,\forall t\in [0,+\infty).
\end{array}$$
Namely
$$\forall \omega \in R, \mbox{Re}
[\displaystyle\frac {b(j\omega )}{a(j\omega )}]>0.$$
This completes the proof.

{\bf  Lemma 7}\ \
Suppose $ a(s)=s^4+a_1s^3+a_2s^2+a_3s+a_4\in H^4,
b(s)=s^4+b_1s^3+b_2s^2+b_3s+a_4\in H^4,$
if $\lambda b(s)+(1-\lambda )a(s)\in H^4,\lambda \in [0,1],$
then $ \Omega _e^a \cap \Omega _e^b\neq \phi $.

{\bf  Proof: }\ \  If
$ \forall \lambda \in [0,1],$
$\lambda b(s)+(1-\lambda )a(s)\in H^4,$
by Lemma 5, for any $ \lambda \in [0,1],$
$$
\begin{array}{ll}
\Omega_e^{a_\lambda}:= & \;\{(x,y)|(a_{\lambda 2}^2-4a_{\lambda 4})x^2
+2(2a_{\lambda 3}-a_{\lambda 1}a_{\lambda 2})xy\\
& \;\ \ \ \  +a_{\lambda 1}^2y^2
-2(a_{\lambda 2}a_{\lambda 3}
-2a_{\lambda 1}a_{\lambda 4})x\\
& \;\ \ \ \  -2a_{\lambda 1}a_{\lambda 3}y+a_{\lambda 3}^2<0\}\\
\end{array}$$
is also an  elliptic region in  the first quadrant of the $x$-$y$ plane, where
$a_{\lambda i}:=a_i+\lambda (b_i-a_i),i=1,2,3,4.$
Apparently, when $\lambda $ changes continuously from $0$ to $1$,
$\Omega_e^{a_\lambda}$ will change continuously from
$\Omega_e^{a}$ to $\Omega_e^{b}.$

Now assume $ \Omega _e^a \cap \Omega _e^b= \phi,$ by Lemma 5
(without loss of generality, suppose
$\displaystyle\frac  {b_3}{b_1}>\displaystyle\frac  {a_3}{a_1}$),
$ \exists v \in [\displaystyle\frac  {a_3}{a_1},
\displaystyle\frac  {b_3}{b_1}]$ and $u\neq 0,$ such that the following line $ l$
\begin{equation}
 l:\ \ \  \displaystyle\frac  {x}{u}+\displaystyle\frac  {y}{v}=1
\end{equation}
 is tangent with $ \Omega _e^a $ and $ \Omega _e^b $ simultaneously, and $l $ separates $ \Omega _e^a $ and $ \Omega _e^b$ (i.e.,
$ \Omega _e^a $ and $ \Omega _e^b$ are on different sides of $l $).

Since $ l$ is tangent with $ \Omega _e^a $, consider
\begin{equation}
\left\{
\begin{array}{l}
\displaystyle\frac  {x}{u}+\displaystyle\frac  {y}{v}=1 \\
(a_2^2-4a_4)x^2+2(2a_3-a_1a_2)xy+a_1^2y^2\\
\quad\quad -2(a_2a_3-2a_1a_4)x-2a_1a_3y+a_3^2=0
\end{array}
\right.
\end{equation}
since $ a(s)$ is Hurwitz stable and $u\neq 0,$
by a direct calculation, we know that the necessary and
sufficient condition for $ l$ being tangent with $ \Omega _e^a $ is
\begin{equation}
uv^2-a_1v^2-a_2uv+a_3v+a_4u=0
\end {equation}

Since $ l$ is tangent with $ \Omega _e^b $, for the same reason, we have
\begin{equation}
uv^2-b_1v^2-b_2uv+b_3v+b_4u=0
\end {equation}
From $ (3)$  and $ (4),$ we obviously have $\forall \lambda \in [0,1],$
\begin{equation}
uv^2-a_{\lambda 1}v^2-a_{\lambda 2}uv+a_{\lambda 3}v+a_{\lambda 4}u=0
\end {equation}
$(5)$ shows that $ l$ is also tangent with $ \Omega _e^{a_{\lambda }}(
\forall \lambda \in [0,1])$, but
$l $ separates $ \Omega _e^a $ and $ \Omega _e^b, $ and when
$\lambda $ changes continuously from $0$ to $1$, $\Omega_e^{a_\lambda}$ will change continuously from $\Omega_e^{a}$ to $\Omega_e^{b},$ which is obviously impossible.This completes the proof.

{\bf Lemma 8}\ \
If $F=\{a_i(s),i=1,2,3,4.\} $
is the set of the four Kharitonov vertex polynomials of a fourth order
 stable interval polynomial family, then
$ \Omega ^{a_2} \subset \Omega ^{a_4} $
and $\Omega ^{a_3} \subset \Omega ^{a_1} .$

{\bf Proof: }\ \
By the definition of the notation $\Omega ^{a},$ it is easy to see that
$$
\begin{array}{ll}
\Omega^{a_1}=\{(x,y)| & \;(a_1^{+}-x)t^2+(a_2^{+}x-a_1^{+}y-a_3^{-})t\\
& \;\quad +(a_3^{-}y-a_4^{-}x)> 0,\forall t\in [0,\infty)\}\\
\Omega^{a_2}=\{(x,y)| & \;(a_1^{-}-x)t^2+(a_2^{-}x-a_1^{-}y-a_3^{+})t\\
& \;\quad +(a_3^{+}y-a_4^{+}x)> 0,\forall t\in [0,\infty)\}\\
\Omega^{a_3}=\{(x,y)| & \;(a_1^{+}-x)t^2+(a_2^{-}x-a_1^{+}y-a_3^{-})t\\
& \;\quad +(a_3^{-}y-a_4^{+}x)> 0,\forall t\in [0,\infty)\}\\
\Omega^{a_4}=\{(x,y)| & \;(a_1^{-}-x)t^2+(a_2^{+}x-a_1^{-}y-a_3^{+})t\\
& \;\quad +(a_3^{+}y-a_4^{-}x)> 0,\forall t\in [0,\infty)\}\\
\end{array}
$$
Obviously, we have $ \Omega ^{a_2} \subset \Omega ^{a_4} $
and $\Omega ^{a_3} \subset \Omega ^{a_1} .$ This completes the proof.

{\bf Lemma 9}\ \
If $F=\{a_i(s),i=1,2,3,4.\} $
is the set of the four Kharitonov vertex polynomials of a fourth order
 stable interval polynomial family, then
$ \cap _{i=1}^4 \Omega ^{a_i} \neq \phi .$

Lemma 9 plays an important role in proving Anderson's claim on
robust SPR synthesis for the fourth-order stable interval
polynomials. For a complete understanding of it, we give three
different proofs in the sequel.

{\bf Proof 1: }\ \
By Lemma 8, we only need to prove that $\Omega ^{a_2} \cap \Omega ^{a_3}\neq \phi .$
By Lemma 7, we know that $\Omega _e^{a_2} \cap \Omega _e^{a_3}\neq \phi ,$
but $\Omega ^{a_2} =\Omega _e^{a_2}\cup \Omega _t^{a_2}$ and
$\Omega ^{a_3} =\Omega _e^{a_3}\cup \Omega _t^{a_3}, $ thus
$\Omega ^{a_2} \cap \Omega ^{a_3}\neq \phi .$ This completes the proof.

{\bf Proof 2: }\ \
Since $F$ is the set of the four Kharitonov vertex
polynomials of a fourth order
 stable interval polynomial family, by Lemma 5,
in the $x$-$y$ plane,
$\Omega_e^{a_2}$ and $\Omega_e^{a_4}$ are both tangent with $x=0$ at
$(0, \displaystyle\frac  {a_3^{+}}{a_1^{-}})$ (denote this tangent point as $A_{24}$);
$\Omega_e^{a_1}$ and $\Omega_e^{a_3}$ are both tangent with $x=0$ at
$(0, \displaystyle\frac  {a_3^{-}}{a_1^{+}})$ (denote this tangent point as $A_{13}$). Denote the tangent point of $\Omega_e^{a_2}$
( $\Omega_e^{a_4}$ ) and $x=a_1^{-}$
as $A_2(a_1^{-},a_2^{-}-\displaystyle\frac  {a_3^{+}}{a_1^{-}})$ (
$A_4(a_1^{-},a_2^{+}-\displaystyle\frac  {a_3^{+}}{a_1^{-}})$ );
and denote the tangent point of $\Omega_e^{a_1}$ ( $\Omega_e^{a_3}$ ) and $x=a_1^{+}$
as $A_1(a_1^{+},a_2^{+}-\displaystyle\frac  {a_3^{-}}{a_1^{+}})$ (
$A_3(a_1^{+},a_2^{-}-\displaystyle\frac  {a_3^{-}}{a_1^{+}})$ ).
Furthermore, denote the intersection points of $x=a_1^{-}$ and the straight line $a_3^{+}y-a_4^{+}x=0,a_3^{+}y-a_4^{-}x=0$ as
$B_2(a_1^{-},\displaystyle\frac  {a_1^{-}a_4^{+}}{a_3^{+}}),
B_4(a_1^{-},\displaystyle\frac  {a_1^{-}a_4^{-}}{a_3^{+}}),$ respectively;
and denote the intersection points of $x=a_1^{+}$ and the straight lines
$a_3^{-}y-a_4^{-}x=0,
a_3^{-}y-a_4^{+}x=0$ as
$B_1(a_1^{+},\displaystyle\frac  {a_1^{+}a_4^{-}}{a_3^{-}}),
B_3(a_1^{+},\displaystyle\frac{a_1^{+}a_4^{+}}{a_3^{-}}),$ respectively.

In what follows, $(A,B)$ stands for the set of points
in the line segment connecting the point $A$ and
the point $B$, not including the endpoints $A$ and $B,$
$[A,B)$ stands for  the set of points in the line segment connecting the point $A$ and
the point $B$, including the endpoint $A$, but not $B,$
$(A,B]$ stands for  the set of points in the line segment connecting the point $A$ and
the point $B$, including the endpoint $B$, but not $A$. Then it is easy to see that
$[A_2,B_2) \subset \Omega ^{a_2},
[A_2,B_2)\subset [A_4,B_4) \subset \Omega ^{a_4},
[A_3,B_3) \subset \Omega ^{a_3},
[A_3,B_3)\subset [A_1,B_1) \subset \Omega ^{a_1},$ and
$(A_{24},A_2]\subset \Omega ^{a_2},(A_{24},A_2]\subset \Omega ^{a_4},$
$(A_{13},A_3]\subset \Omega ^{a_3},(A_{13},A_3]\subset \Omega ^{a_1}.$

Denote $A_3^{\star}$ as $(a_1^{-},(\displaystyle\frac  {a_2^{-}}{a_1^{+}}-2
\displaystyle\frac  {a_3^{-}}{{a_1^{+}}^2})a_1^{-}+
\displaystyle\frac  {a_3^{-}}{a_1^{+}}),$ then $A_3^{\star}\in (A_{13},A_3].$

If $\displaystyle\frac  {a_3^{+}}{a_1^{-}}=
\displaystyle\frac  {a_3^{-}}{a_1^{+}},$ i.e., $ a_1^{-}=a_1^{+}$ and
$ a_3^{-}=a_3^{+}.$ Then, take $\delta >0,\delta$ sufficiently small, by Lemma 5, it is easy to
verify that $(\delta, \displaystyle\frac  {a_3^{+}}{a_1^{-}})
\in \cap _{i=1}^4 \Omega _e^{a_i},$ thus
$ \cap _{i=1}^4 \Omega ^{a_i} \neq \phi .$

Now, suppose $\displaystyle\frac  {a_3^{+}}{a_1^{-}}>
\displaystyle\frac  {a_3^{-}}{a_1^{+}}$
and
$$a_2^{-}-\displaystyle\frac  {a_3^{+}}{a_1^{-}}
\geq (\displaystyle\frac  {a_2^{-}}{a_1^{+}}-2
\displaystyle\frac  {a_3^{-}}{{a_1^{+}}^2})a_1^{-}+
\displaystyle\frac  {a_3^{-}}{a_1^{+}}$$

\noindent It is easy to verify that
$$(\displaystyle\frac  {a_2^{-}}{a_1^{+}}-2
\displaystyle\frac  {a_3^{-}}{{a_1^{+}}^2})a_1^{-}+
\displaystyle\frac  {a_3^{-}}{a_1^{+}}>
\displaystyle\frac  {a_1^{-}a_4^{+}}{a_3^{+}}$$
Thus, we have $A_3^{\star}\in [A_2,B_2). $ Hence $ A_3^{\star}\in [A_2,B_2)\cap (A_{13},A_3].$ Therefore
$A_3^{\star}\in \cap _{i=1}^4 \Omega ^{a_i}.$ Thus
$ \cap _{i=1}^4 \Omega ^{a_i} \neq \phi .$

\noindent Finally, with $\displaystyle\frac  {a_3^{+}}{a_1^{-}}>
\displaystyle\frac  {a_3^{-}}{a_1^{+}},$
if
$$a_2^{-}-\displaystyle\frac  {a_3^{+}}{a_1^{-}}
< (\displaystyle\frac  {a_2^{-}}{a_1^{+}}-2
\displaystyle\frac  {a_3^{-}}{{a_1^{+}}^2})a_1^{-}+
\displaystyle\frac  {a_3^{-}}{a_1^{+}}$$
then it is easy to see that $(A_{13},A_3]\cap (A_{24},A_2]\neq \phi$
and $(A_{13},A_3]\cap (A_{24},A_2]\subset \cap _{i=1}^4 \Omega ^{a_i}.$
Thus, we also have $ \cap _{i=1}^4 \Omega ^{a_i} \neq \phi .$
This completes the proof.

{\bf Proof 3: }\ \
$A_{13}(0, \displaystyle\frac  {a_3^{-}}{a_1^{+}}),
A_{24}(0, \displaystyle\frac  {a_3^{+}}{a_1^{-}}),
B_2(a_1^{-},\displaystyle\frac  {a_1^{-}a_4^{+}}{a_3^{+}})$ and $B_3
(a_1^{+},\displaystyle\frac{a_1^{+}a_4^{+}}{a_3^{-}})$
are defined identically as in the Proof 2 above.
$(A,B)$ stands for the set of points
in the line segment connecting the point $A$ and
the point $B$, but not including the endpoints $A$ and $B.$

If $\displaystyle\frac  {a_3^{+}}{a_1^{-}}=
\displaystyle\frac  {a_3^{-}}{a_1^{+}},$ i.e., $ a_1^{-}=a_1^{+}$ and
$ a_3^{-}=a_3^{+}.$ Then, take $\delta >0,\delta$ sufficiently small,
by Lemma 5, it is easy to
verify that $(\delta, \displaystyle\frac  {a_3^{+}}{a_1^{-}})
\in \cap _{i=1}^4 \Omega _e^{a_i},$ thus
$ \cap _{i=1}^4 \Omega ^{a_i} \neq \phi .$

Now, suppose $\displaystyle\frac  {a_3^{+}}{a_1^{-}}>
\displaystyle\frac  {a_3^{-}}{a_1^{+}},$
then it is easy to see that $(A_{13},B_3)\cap (A_{24},B_2)\neq \phi$
and $(A_{13},B_3)\cap (A_{24},B_2)\subset \cap _{i=1}^4 \Omega ^{a_i}.$
Thus, we also have $ \cap _{i=1}^4 \Omega ^{a_i} \neq \phi .$
This completes the proof.

{\bf Lemma {\bf 10}}\ \
Suppose $ a(s)=s^4+a_1s^3+a_2s^2+a_3s+a_4\in H^4,
b(s)=s^3+xs^2+ys+z, $ and $\forall \omega \in R,
\mbox{Re}
[\displaystyle\frac {b(j\omega )}{a(j\omega )}]>0,$  take
\begin{center}
$ \stackrel{\sim }{b}(s):=b(s)+r\cdot c(s),\ \ r>0,r$ sufficiently small
\end{center}
where $c(s)$ is a fixed fourth-order monic polynomial. Then
$ \displaystyle\frac {\stackrel{\sim }{b}(s)}{a(s)}$
is strictly positive real.

{\bf Proof: }\ \ Obviously, $ \partial (a)=
\partial ({\stackrel{\sim }{b}}),$ namely,  $ {\stackrel{\sim }{b}}(s)$ and $ a(s)$ have the same order. Since $ a(s)\in H^4,$
there exists $ \omega _1 >0$ such that, for all
$ \mid \omega \mid \geq \omega _1,$
Re$ [\displaystyle\frac {{\stackrel{\sim }{b}}(j\omega )}{a(j\omega )}]>0.$
Denote
\begin{center}
$ M_1=\dinf_{\mid \omega \mid \leq \omega _1}\mbox{Re} [ \displaystyle\frac
{b(j\omega )}{a(j\omega )}] $ \quad \quad
$ N_1=\dsup_{\mid \omega \mid \leq \omega _1}\displaystyle {\mid
\mbox{Re}[\displaystyle\frac
{c(j\omega )}{a(j\omega )}] \mid }$
\end{center}
Then $  M_1>0 $  and $  N_1>0.$  Choosing
$  0<r<\displaystyle\frac {M_1}{N_1},$ it can be directly verified that
$$ \forall \omega \in R, \mbox{Re}
[\displaystyle\frac {{\stackrel{\sim }{b}}(j\omega )}{a(j\omega )}]>0$$
This completes the proof.

Now Theorem 3 is proved by simply combining Lemmas 5-7 and Lemma 10.
Theorem 4 is proved by simply combining Lemmas 5-6 and Lemmas 9-10.
\begin{center}
\section*{4. DISCUSSIONS AND EXAMPLES}
\end{center}
The following three examples correspond to different cases in the proof of our main results.

{\bf Example 1} \ \
Suppose
$ a_1(s)=s^4+89s^3+56s^2+88s+1,
 a_2(s)=s^4+11s^3+56s^2+88s+50,
 a_3(s)=s^4+89s^3+56s^2+88s+50,
 a_4(s)=s^4+11s^3+56s^2+88s+1$
are the four Kharitonov vertex polynomials of a fourth-order interval polynomial set $K$, it is easy to check
using Kharitonov's Theorem that $K$ is robustly stable.
By our method as in the constructive proof of
Theorem 4, it is easy to get $(11,7.6657)\in \cap _{i=1}^4 \Omega ^{a_i}.$
Thus, choose $b(s)= s^3+11s^2+7.76657s+\varepsilon,$
where $\varepsilon $ is a sufficiently small positive number
($\varepsilon $ is determined by Lemma 6, in this example,
$0<\varepsilon \leq 3$), take $\varepsilon =2,$ by Lemma 6,
$\forall \omega \in R,
\mbox{Re}
[\displaystyle\frac {b(j\omega )}{a_i(j\omega )}]>0,i=1,2,3,4.$
Finally, let $\stackrel{\sim }{b}(s):=b(s)+r\cdot s^4,$
where $r>0,r$ sufficiently small ($r$ is determined by Lemma 10, in this example,
$0<r \leq 0.5$),
it is easy to check that $ \displaystyle\frac {\stackrel{\sim }{b}(s)}{a_i(s)},
i=1,2,3,4,$
are strictly positive real
(note that $b(s)$ and ${\stackrel{\sim }{b}(s)}$ are not unique).

In this example, if we take $F=\{a_1(s)=s^4+11s^3+56s^2+88s+50,
a_2(s)=s^4+89s^3+56s^2+88s+50\},$
then it is exactly the counter-example provided in [17].
It can be checked that $F$ does not satisfy
the sufficient conditions in [10,16,17], but we can use the methods in
[13-15,19,20] to do SPR synthesis.
When $F$ is enlarged to the interval polynomial set  $K$ in
this example, the synthesis methods in [13-15] fail too,
but we can still use the methods in [19,20]
to do synthesis. It is quite straightforward to do
synthesis using the method in this paper.

{\bf Example 2} \ \
Suppose
$ a_1(s)=s^4+5s^3+6s^2+4s+0.5,
 a_2(s)=s^4+2s^3+6s^2+6s+1,
 a_3(s)=s^4+5s^3+6s^2+4s+1,
 a_4(s)=s^4+2s^3+6s^2+6s+0.5$
are the four Kharitonov vertex polynomials of a fourth-order interval polynomial set $K$, it is easy to check
using Kharitonov's Theorem that $K$ is robustly stable.
By our method as in the constructive proof of
Theorem 4, it is easy to get $(2,2.56)\in \cap _{i=1}^4 \Omega ^{a_i}.$
Thus, choose $b(s)= s^3+2s^2+2.56s+\varepsilon,$
where $\varepsilon $ is a sufficiently small positive number
(in this example,
$0<\varepsilon \leq 1$), take $\varepsilon =0.5,$ by Lemma 6,
$\forall \omega \in R,
\mbox{Re}
[\displaystyle\frac {b(j\omega )}{a_i(j\omega )}]>0,i=1,2,3,4.$
Finally, let $\stackrel{\sim }{b}(s):=b(s)+r\cdot s^4,$
where $r>0,r$ sufficiently small (in this example,
$0<r \leq 0.5$),
it is easy to check that $ \displaystyle\frac {\stackrel{\sim }{b}(s)}{a_i(s)},
i=1,2,3,4,$
are strictly positive real.

{\bf Example 3} \ \
Suppose
$ a_1(s)=s^4+2.5s^3+6s^2+4s+0.5,
 a_2(s)=s^4+2s^3+5s^2+6s+5,
 a_3(s)=s^4+2.5s^3+5s^2+4s+5,
 a_4(s)=s^4+2s^3+6s^2+6s+0.5$
are the four Kharitonov vertex polynomials of a fourth-order interval polynomial set $K$, it is easy to check
using Kharitonov's Theorem that $K$ is robustly stable.
By our method as in the constructive proof of
Theorem 4, it is easy to get $(1.1475,2.4262)
\in \cap _{i=1}^4 \Omega ^{a_i}.$
Thus, choose $b(s)= s^3+1.1475s^2+2.4262s+\varepsilon,$
where $\varepsilon $ is a sufficiently small positive number
(in this example,
$0<\varepsilon \leq 1$), take $\varepsilon =0.5,$ by Lemma 6,
$\forall \omega \in R,
\mbox{Re}
[\displaystyle\frac {b(j\omega )}{a_i(j\omega )}]>0,i=1,2,3,4.$
Finally, let $\stackrel{\sim }{b}(s):=b(s)+r\cdot s^4,$
where $r>0,r$ sufficiently small (in this example,
$0<r \leq 0.2$),
it is easy to check that $ \displaystyle\frac {\stackrel{\sim }{b}(s)}{a_i(s)},
i=1,2,3,4,$
are strictly positive real.

{\bf Remark  {\bf 1}} \ \
From the proofs of Theorem 3 and Theorem 4, we can see that,
this paper not only proves the existence,
but also provides a design procedure.

{\bf Remark  {\bf 2}} \ \
Lemma 10 actually holds for arbitrary $n$-th order polynomials$^{[19,20]}$.%

{\bf Remark  {\bf 3}} \ \
The constructive synthesis method is also insightful and
helpful in solving the general robust SPR synthesis
problem. In fact, we have recently succeeded in proving the existence on  robust SPR synthesis
for fifth-order stable convex combinations
using a similar method$^{[19]}$.
The SPR synthesis for higher-order systems is currently under investigation.

{\bf Remark {\bf 4}} \ \
Robust stability of a polynomial segment can be checked by many efficient methods, e.g., eigenvalue method, root locus method,
value set method, etc.$^{[8,9]}$.
Robust stability of $K$ in Theorem 4 can be ascertained by
checking only two Kharitonov vertex polynomials$^{[24]}$.

{\bf Remark  {\bf 5}} \ \
From the proofs of Lemma 9, we can establish the relationship
between SPR synthesis for the fourth-order polynomial
segments and SPR synthesis for the fourth-order interval polynomials.
In fact, it is easy to see that Theorem 3 implies Theorem 4.
Similarly, Theorem 1 implies Theorem 2.
However, similar results may not be true for higher-order ($n\geq 5$) systems.
This subject is currently under investigation.

{\bf  Remark {\bf 6}} \ \
Our results can easily be generalized to discrete-time case. %

Finally, it should also be pointed out that, for the vertex set
$F=\{a_i(s)=s^n+\sum_{l=1}^na_l^{(i)}s^{n-l},i=1,2,\cdots,m.\}$
of a general polytopic polynomial family $\overline{F}$, even if
$\overline{F}$ is robustly stable, it is still possible that there does not exist a polynomial $c(s)\in H^{n-1},$ such that,
$\forall \omega \in R,
\mbox{Re}
[\displaystyle\frac {c(j\omega )}{a(j\omega )}]>0,
$ for all $a(s)\in \overline{F}$.

To see this, let us look at an example of a third order triangle polynomial family.

{\bf  Example  4} \ \
Let
$F=\{ a_1(s)=s^3+2.6s^2+37s+64,$
$a_2(s)=s^3+17s^2+83s+978,$
$a_3(s)=s^3+15s^2+28s+415 \}.$
It is easy to verify that
$a_i(s),i=1,2,3,$ are Hurwitz stable. Moreover, all edges of $\overline{F}$, i.e.,
$\lambda a_i(s)+(1-\lambda)a_j(s), \lambda \in [0,1],i,j=1,2,3,$ are also Hurwitz stable. Therefore, by Edge Theorem$^{[6-9]}$,
$\overline{F}$ is robustly stable.
On the other hand, by a direct computation,  we can easily see that
there does not exist a polynomial $c(s)\in H^2,$
such that
$\forall \omega \in R,
\mbox{Re}
[\displaystyle\frac {c(j\omega )}{a_i(j\omega )}]>0,
i=1,2,3.$

Note that, in this example, though
there does not exist a polynomial $c(s)\in H^2$
such that
$\forall \omega \in R,
\mbox{Re}
[\displaystyle\frac {c(j\omega )}{a_i(j\omega )}]>0,
i=1,2,3.$
But if we take $\stackrel{\sim }{c}(s)=s^3+6s^2+73s+68,$
it is easy to check
 $ \displaystyle\frac {\stackrel{\sim }{c}(s)}{a_i(s)},
i=1,2,3,$
are strictly positive real. This shows some intrinsic differences between the SPR synthesis
of interval polynomial families and the SPR synthesis of polytopic polynomial families. This problem
deserves further investigation.
\vskip 7mm
\begin{center}
\section*{5. CONCLUSIONS}
\end{center}
\vskip -8mm {\small We have proved that, for low-order ($n\leq 4$)
stable polynomial segments or interval polynomials, there always
exists a fixed polynomial such that their ratio is SPR-invariant,
thereby providing a rigorous proof of Anderson's claim on SPR
synthesis for the fourth-order stable interval polynomials.
Moreover, the relationship between SPR synthesis for low-order
polynomial segments and SPR synthesis for low-order interval
polynomials has also been discussed.}
\subsection*{\centering Acknowledgments}
\vskip -0.1in
\noindent {\small This work was supported by the National Key Project of China,
the National Natural Science Foundation of China (69925307),
Natural Science Foundation of Chinese Academy of Sciences and
National Lab of
Intelligent Control and Systems of Tsinghua University.}

\end{document}